\documentclass{amsart}

\newtheorem{theorem}{Theorem}[section]

\newtheorem{corollary}[theorem]{Corollary}
\newtheorem{lemma}[theorem]{Lemma}

\theoremstyle{remark}
\newtheorem*{remark}{Remark}
\theoremstyle{definition}

\newtheorem{example}[theorem]{Example}
\newcommand{\X}{\mathbf{X}}

\newcommand{\ff}{{\mathbb{ F\!}}}

\newcommand{\x}{{\bf x}}

\usepackage{graphicx}
\usepackage{amsfonts}
\begin{document}
\numberwithin{equation}{section}

% Title
\title[Linear recurrences of symmetric boolean functions]{Linear recurrences and asymptotic behavior of exponential sums of symmetric boolean functions}
 
% Information about authors
\author{Francis N. Castro}
\address{Department of Mathematics, University of Puerto Rico, San Juan, PR 00931}
\email{franciscastr@gmail.com}

\author{Luis A. Medina}
\address{Department of Mathematics, University of Puerto Rico, San Juan, PR 00931}
\email{luis.medina@uprrp.edu}

% Abstract
\begin{abstract}
In this paper we give an improvement of the degree of the homogeneous linear recurrence with integer coefficients that exponential sums  of symmetric Boolean functions satisfy. This improvement is tight.  We also compute the asymptotic behavior of symmetric Boolean functions and provide a formula that allows us to determine if a symmetric boolean function is asymptotically not balanced.  In particular, when the degree of the symmetric function is a power of two, then the exponential sum is much smaller than $2^n$.
\end{abstract}

% General info
\date{\today}
% Keywords

\maketitle

%%%%%%%%%%%%%%%%%%%%%%%%%%%%%%%%%%%%%%%%%%%%
% Section: Introduction
%%%%%%%%%%%%%%%%%%%%%%%%%%%%%%%%%%%%%%%%%%%%
\section{Introduction}
Boolean functions are one of the most studied objects in mathematics. They are important in many applications, for example, in the design of stream ciphers, block  and hash functions. These functions also play a vital role in cryptography as they are used as filter and combination generator of stream ciphers based on linear feed-back shift registers. The case of boolean functions of degree 2 has been intensively studied because of its relation to bent functions (for example see \cite{Roth}, \cite{Bey}).

One can find many papers and books discussing the properties of boolean functions (see for  examples \cite{carlet}, \cite{cusick}, \cite{canf} and \cite{carlet1}). The  subject can be studied from the point of view of complexity theory or from the algebraic point of view as we do in this paper, where we compute the asymptotic behavior of exponential sums of  symmetric boolean functions.

The correlation between two Boolean functions of $n$ inputs is defined as the number of times the functions agree minus the number of times they disagree all divided by $2^n$, i.e., 
\begin{equation}
C(F_1,F_2)=\frac{1}{2^n}\sum_{x_1,\ldots,x_n\in \{0,1\}} (-1)^{F_1(x_1,\ldots,x_n)+F_2(x_1,\ldots,x_n)}.
\end{equation}
In this paper we are interested in the case when $F_1$ and $F_2$ are symmetric boolean functions. Without loss of generality, we write $C(F)$ instead of $C(F_1,F_2)$, where $F$ is a symmetric boolean function. In \cite{can1}, A. Canteaut and M. Videau studied in details symmetric boolean functions. They established a link between the periodicity of the simplified value vector of a symmetric Boolean function and its degree. %Their result reduces the amount of memory required for representing a symmetric function and has some consequences from a cryptographic point of view. 
They also determined all balanced symmetric functions of degree less than or equal to 7. In \cite{von}, J. von zur Gathen and J. Rouche  found  all the balanced symmetric boolean functions up to 128 variables. 

In\cite{cai}, J. Cai et. al. computed a closed formula for the correlation between any two symmetric Boolean functions. This formula implies that $C(F)$ satisfies a homogeneous linear recurrence with integer coefficients and provides an upper bound for the degree of the minimal recurrence of this type that $C(F)$ satisfies. In this paper we give an improvement to the degree of the minimal homogeneous linear recurrence with integer coefficients satisfying by $C(F)$. In particular, our lower and upper bounds are tight in many cases. Also, in the case of an elementary symmetric function we provide the minimal linear recurrence it satisfies.  

We  also compute the asymptotic  value of $C(F)$. In particular, we give infinite families of boolean functions that are asymptotically not balanced, i.e., $\lim_{n\rightarrow \infty} C(F)\not=0$. In the case that $F$ is an elementary symmetric function, $F$ is asymptotically not balanced if and only if its degree is not a power of 2.   When the asymptotic value of $C(F)$ is zero, we compute the asymptotic values of
\begin{equation}
\frac{1}{|\lambda|^n} \sum_{x_1,\ldots,x_n\in \{0,1\}} (-1)^{F(x_1,\ldots,x_n)},
\end{equation}
where  $\lambda$ and $\bar{\lambda}$ are the roots with the biggest modulus of the characteristic polynomial associated to the exponential sum of $F$.  We prove that the coefficient of $\lambda$ is not identically zero and obtain information about the spectrum of $F(X_1,\ldots,X_n)$. In particular, its limit  is a periodic function in $n$.

%%%%%%%%%%%%%%%%%%%%%%%%%%%%%%%%%%%%%%%%%%%%%%%%%%%%%%%%%%%%%%%%%
% Preliminaries
%%%%%%%%%%%%%%%%%%%%%%%%%%%%%%%%%%%%%%%%%%%%%%%%%%%%%%%%%%%%%%%%%
\section{Preliminaries}
Let $\ff$ be the binary field, $\ff\,^n=\{(x_1,\ldots,x_n)\,|\, x_i\in \ff, i=1,\ldots,n\}$, and $F(\X)=F(X_1,\ldots,X_n)$ be a polynomial in $n$ variables over $\ff$. The exponential sum associated to $F$ over $\ff$ is:
\begin{equation}\label{eq1}
 S(F)=\sum_{\x\in \ff^{\,n}}(-1)^{F(\x)}.
\end{equation}
Note that $S(F)=2^nC(F)$.  A boolean function $F({\bf X})$ is called balanced if $S(F)=0$.  This property is important for some applications in cryptography.  P. Sarkar and S. Moitra \cite{sarkar} found a lower bound for the number of symmetric balanced boolean functions.  This number is bigger than or equal to $2^{(n+1)/2}+2^{(n+1)/2-3}$ where $n\geq3$ is the number of variables.

In this paper we study exponential sums associated to symmetric boolean functions $F$.  Any symmetric function is a linear combination of elementary symmetric polynomials, thus we start with exponential sums of elementary symmetric polynomials.

Let $\sigma_{n,k}$ be the elementary symmetric polynomial in $n$ variables of degree $k$. For example,
\begin{equation}
\sigma_{4,3} = X_1 X_2 X_3+X_1 X_4 X_3+X_2 X_4 X_3+X_1 X_2 X_4.
\end{equation}
Fix $k\geq 2$ and let $n$ vary. Consider the sequence of exponential sums $\{S(\sigma_{n,k})\}_{n\in \mathbb{N}}$ where
\begin{equation}
\label{expsumsymm}
S(\sigma_{n,k})=\sum_{x_1,\cdots,x_n\in \mathbb{F}} (-1)^{\sigma_{n,k}(x_1,\cdots,x_n)}.
\end{equation}
Define $A_j$ to be the set of all $(x_1,\cdots,x_n)\in \mathbb{F}^{n}$ with exactly $j$ entries equal to 1.  Clearly, $|A_j|=\binom{n}{j}$ and $\sigma_{n,k}({\bf x})=\binom{j}{k}$ for ${\bf x}\in A_j$.  Therefore,
\begin{equation}
\label{main}
S(\sigma_{n,k})=\sum_{j=0}^n (-1)^{\binom{j}{k}}\binom{n}{j}.
\end{equation}
In general, if $1\leq k_1 < k_2 < \cdots < k_s$ are fixed integers, then
\begin{equation}
\label{maingen}
S(\sigma_{n,k_1}+\sigma_{n,k_2}+\cdots +\sigma_{n,k_s}) =\sum_{i=0}^n (-1)^{\binom{j}{k_1}+\binom{j}{k_2}+\cdots+\binom{j}{k_s}}\binom{n}{j}.
\end{equation}

\begin{remark} 
Note that the sum on the right hand side of (\ref{maingen}) makes sense for values of $n$ less than $k_s$, while $S(\sigma_{n,k_1}+\cdots+\sigma_{n,k_s})$ does not.  However, throughout the paper we let $S(\sigma_{n,k_1}+\cdots+\sigma_{n,k_s})$ to be defined by the sum in (\ref{maingen}), even for values of $n$ less than $k_s$.
\end{remark}

%%%%%%%%%%%%%%%%%%%%%%%%%%%%%%%%%%%%%%%%%%%%
% Section: The recurrence
%%%%%%%%%%%%%%%%%%%%%%%%%%%%%%%%%%%%%%%%%%%%
\section{The recurrence}
Computer experimentation suggests that for fix $1\leq k_1 < \cdots <k_s$, the sequence $\{S(\sigma_{n,k_1}+\cdots+\sigma_{n,k_s})\}_{n\in \mathbb{N}}$ satisfies a homogeneous linear recurrence with integer coefficients.  For example, if we consider $\{S(\sigma_{n,7})\}_{n\in \mathbb{N}}$ and type 
\begin{eqnarray*}
\texttt{FindLinearRecurrence[Table[Sum[((-1)\^{}Binomial[m,7])*}\\\texttt{Binomial[n,m],\{m,0,n\}],\{n,1,30\}]]}
\end{eqnarray*}
into {\it Mathematica} 7, then it returns
\begin{eqnarray*}
\texttt{\{8,-28,56,-70,56,-28,8\}}.
\end{eqnarray*}
This suggests that $\{S(\sigma_{n,7})\}_{n\in \mathbb{N}}$ satisfies the recurrence
\begin{equation}
\label{firstex}
x_n = 8x_{n-1}-28x_{n-2}+56x_{n-3}-70x_{n-4}+56x_{n-5}-28x_{n-6}+8x_{n-7}.
\end{equation}
If we continue with these experiments, we arrive to the observation that if $r=\lfloor\log_2(k_s)\rfloor+1$, then $\{S(\sigma_{n,k_1}+\cdots+\sigma_{n,k_s})\}_{n\in \mathbb{N}}$ seems to satisfy the recurrence
\begin{equation}
\label{altrec}
x_{n}=\sum_{m=1}^{2^r-1}(-1)^{m-1}\binom{2^r}{m}x_{n-m}.
\end{equation}

This result can be proved using elementary machinery.  The idea is to use the fact that if $r=\lfloor \log_2(k_s) \rfloor+1$, then
\begin{equation}
\label{frombino}
\binom{j+i2^{r}}{k_m} \equiv \binom{j}{k_m} (\text{mod }2)
\end{equation}
for all non-negative integer $i$ and $m=1,2,\cdots, s,\,$ to show inductively that the family of sequences
\begin{equation}
\label{family}
a_{n,r,i}=\sum_{j}\binom{n}{2^rj+i}=\sum_{j \equiv i \,\, (\text{mod }2^r)}\binom{n}{j},
\end{equation}
$i=0,1,\cdots 2^r-1$ satisfies the same recurrence (\ref{altrec}).  However, we should point out the fact that $\{S(\sigma_{n,k_1}+\cdots+\sigma_{n,k_s})\}_{n\in \mathbb{N}}$ satisfying (\ref{altrec}) is a consequence of the following theorem of J. Cai et al \cite{cai}. 
\begin{theorem}
\label{invariant}
Fix $1\leq k_1< \cdots <k_s$ and let $r=\lfloor \log_2(k_s) \rfloor+1$. The value of the exponential sum $S(\sigma_{n,k_1}+\cdots+\sigma_{n,k_s})$ is given by
\begin{eqnarray}
\label{genvalue}\nonumber
S(\sigma_{n,k_1}+\cdots+\sigma_{n,k_s}) &=& \sum_{i=0}^n (-1)^{\binom{i}{k_1}+\cdots+\binom{i}{k_s}}\binom{n}{i} \\
&=& c_0(k_1,\cdots,k_s) 2^n + \sum_{j=1}^{2^r-1} c_j(k_1,\cdots,k_s) (1+\zeta_j)^n,
\end{eqnarray}
where $\zeta_j = \exp\left(\frac{\pi \sqrt{-1}\, j}{2^{r-1}}\right)$ and 
\begin{equation}
\label{coeffs}
c_j(k_1,\cdots,k_s)=\frac{1}{2^r}\sum_{i=0}^{2^r-1}(-1)^{\binom{i}{k_1}+\cdots+\binom{i}{k_s}}\zeta_j^{-i}.
\end{equation}
\end{theorem}
The proof of Cai et al. rely on linear algebra.  They wrote
\begin{equation}
S(\sigma_{n,k_1}+\cdots+\sigma_{n,k_s}) = \sum_{i=0}^{2^r-1}(-1)^{\binom{i}{k_1}+\cdots+\binom{i}{k_s}}a_{n,r,i},
\end{equation}
and used the elementary identity
\begin{equation}
\binom{n}{k}=\binom{n-1}{k}+\binom{n-1}{k-1},
\end{equation}
to find a recurrence for
\begin{equation}
{\bf a}_{n,r}=\left(
\begin{array}{c}
a_{n,r,1} \\
a_{n,r,2} \\
\vdots\\
a_{n,r,2^{r}-1}
\end{array}
\right)
\end{equation}
of the form ${\bf a}_{n,r} ={\bf M a}_{n-1,r}$, for some matrix ${\bf M}$.  Finding the eigenvalues and corresponding eigenvectors of {\bf M}, they were able to solve this recurrence and prove Theorem \ref{invariant}.  See \cite{cai} for more details.

From Theorem \ref{invariant} it is now evident that $\{S(\sigma_{n,k_1}+\cdots+\sigma_{n,k_s})\}_{n\in \mathbb{N}}$ satisfies (\ref{altrec}).  Moreover, the roots of the characteristic polynomial associated to the linear recurrence (\ref{altrec}) are all different and the polynomial is given by
\begin{eqnarray}
P_r(x)&=&\sum_{m=0}^{2^r-1}(-1)^{m}\binom{2^r}{m}x^{2^r-1-m}\\ \nonumber
&=& (x-2)\Phi_4(x-1)\Phi_8(x-1)\cdots\Phi_{2^r}(x-1),
\end{eqnarray}
where $\Phi_m(x)$ represents the $m$-th cyclotomic polynomial
\begin{equation}
\Phi_m(x)=\prod_{\zeta^m=1\text{ primitive }} (x-\zeta).
\end{equation}
Even though $\{S(\sigma_{n,k_1}+\cdots+\sigma_{n,k_s})\}_{n\in \mathbb{N}}$ satisfies (\ref{altrec}), in many instances (\ref{altrec}) is not the minimal homogenous linear recurrence with integer coefficients that $\{S(\sigma_{n,k_1}+\cdots+\sigma_{n,k_s})\}_{n\in \mathbb{N}}$ satisfies.  For example, $\{S(\sigma_{n,3}+\sigma_{n,5})\}_{n\in \mathbb{N}}$ satisfies (\ref{firstex}), but its minimal recurrence is
\begin{equation}
x_n = 6x_{n-1}-14x_{n-2}+16x_{n-3}-10x_{n-4}+4x_{n-5}.
\end{equation}
In the next section we use Theorem \ref{invariant} to give some improvements on the degree of the minimal linear recurrence associated to $\{S(\sigma_{n,k_1}+\cdots+\sigma_{n,k_s})\}_{n\in \mathbb{N}}$.

%%%%%%%%%%%%%%%%%%%%%%%%%%%%%%%%%%%%%%%%%%%%%%%%%%%%%%%%%%
% Section: On the degree of the recurrence relation
%%%%%%%%%%%%%%%%%%%%%%%%%%%%%%%%%%%%%%%%%%%%%%%%%%%%%%%%%%
\section{On the degree of the recurrence relation}
\label{sectiondegree}
Now that we are equipped with equation (\ref{coeffs}), we move to the problem of reducing the degree of the recurrence relation that our sequences of exponential sums satisfy. The idea behind our approach is very simple.  Consider all roots $1+\zeta$'s of $\Phi_{2^{t+1}}(x-1)$ where $1\leq t\leq r-1$.  We know that $(1+\zeta)^n$ appears in (\ref{genvalue}).  If we show that the coefficient that corresponds to $(1+\zeta)^n$ is zero for each $1+\zeta$, then we reduce the degree of the characteristic polynomial, and therefore the degree of the recurrence, by $2^t$.  

However, note that $\Phi_{2^{t+1}}(x-1)$ is irreducible over $\mathbb{Q}$ (according to Eisenstein's criterion on $\Phi_{2^{t+1}}(x-1)$ with $\Phi_{2^{t+1}}(x)=x^{2^t}+1$, see \cite{garling}).  Therefore, the coefficients related to the roots of $\Phi_{2^{t+1}}(x-1)$ are either all zeros or all non-zeros.  In view of (\ref{coeffs}), this can be determined by checking whether or not the sum
\begin{equation}
\sum_{m=0}^{2^r-1}(-1)^{\binom{m}{k_1}+\cdots+\binom{m}{k_s}}e^{\frac{\pi\sqrt{-1}m}{2^t}}
\end{equation}
is zero.

We discuss first the case of the exponential sum of one elementary symmetric polynomial, i.e. $\{S(\sigma_{n,k}\}_{n\in \mathbb{N}}$.  We start with the following elementary result.
\begin{lemma}[{\bf Lucas' theorem}]
Let $n$ be a natural number with $2$-adic expansion $n=2^{a_1}+2^{a_2}+\cdots+2^{a_l}$.  The binomial coefficient $\binom{n}{k}$ is odd if and only if $k$ is either $0$ or a sum of some of the $2^{a_i}$'s.
\end{lemma}
\begin{proof}
Recall that $(1+x)^{2^m}\equiv 1+x^{2^m}\,(\text{mod }2)$ for all non-negative integer $m$, thus $(1+x)^n\equiv (1+x^{2^{a_1}})(1+x^{2^{a_2}})\cdots(1+x^{2^{a_l}})\,(\text{mod }2)$.  Note that the coefficient of $x^k$ in $(1+x^{2^{a_1}})(1+x^{2^{a_2}})\cdots(1+x^{2^{a_l}})$ is 1 if and only if $k=0$ or a sum of some of the $2^{a_i}$'s.
\end{proof}

\noindent
The next result is an immediate consequence of the above lemma.
\begin{corollary}
\label{modds}
Fix a natural number $k$.  Suppose its $2$-adic expansion is $k=2^{a_1}+2^{a_2}+\cdots+2^{a_l}$.  A natural number $m$ is such that $\binom{m}{k}$ is odd if and only if $m$ has a $2$-adic expansion of the form
\begin{equation}
m=k+\sum_{2^i \not \in \{2^{a_1},2^{a_2},\cdots, 2^{a_l}\}}\delta_i 2^i
\end{equation} 
where $\delta_i \in \{0,1\}$.
\end{corollary}
\begin{remark}
Let $k\geq 1$ be an integer with $2$-adic expansion $k=2^{a_1}+\cdots+2^{a_l}$. Suppose $m \in \{0,1,2,3,\cdots, 2^r-1\}$ is such that $\binom{m}{k}$ is odd. Note that Corollary \ref{modds} implies
\begin{equation}
\label{negativeM}
m=k+\delta_1 2^{b_1}+\delta_2 2^{b_2}+\cdots+\delta_t 2^{b_f},
\end{equation}
where $\{2^{b_1},2^{b_2},\cdots, 2^{b_f}\}=\{1,2,2^2,\cdots, 2^{r-1}\} \backslash\{2^{a_1},2^{a_2},\cdots, 2^{a_l}\}$.
\end{remark}
We now proceed to show which coefficients $c_j(k)$ are zero.  We start with $c_0(k)$.
\begin{lemma}
\label{c0is0}
Suppose $k\geq 2$ is an integer.  Then,
\begin{equation}
\label{valc0}
c_0(k)=\frac{2^{w_2(k)-1}-1}{2^{w_2(k)-1}},
\end{equation}
where $w_2(k)$ is the sum of the binary digits of $k$.  In particular, $c_0(k)=0$ if and only if $k$ is a power of two.
\end{lemma}
\begin{proof}
Recall from Theorem \ref{invariant} that
$$c_0(k)=\frac{1}{2^r}\sum_{m=0}^{2^r-1}(-1)^{\binom{m}{k}},$$
where $r=\lfloor \log_2(k) \rfloor+1$.  Re-write $c_0(k)$ as
\begin{equation}
c_0(k)=\frac{1}{2^r}\left(2^r-2\sum_{m\in N}1\right),
\end{equation}
where 
\begin{equation}
N=\left\{m \in \{0,1,2,3,\cdots, 2^r-1\}: \binom{m}{k} \text{ is odd.}\right\}
\end{equation}
Note that (\ref{negativeM}) implies that the cardinality of $N$ is $2^{r-w_2(k)}$. A simple calculation yields the result.
\end{proof}
\begin{remark}
Lemma \ref{c0is0} shows that if $k$ is not a power of two, then $\sigma_{n,k}$ is not balanced for sufficiently large $n$.  We say that $\sigma_{n,k_1}+\cdots+\sigma_{n,k_s}$ is asymptotically not balanced if
\begin{equation}
\lim_{n\to\infty}\frac{S(\sigma_{n,k_1}+\cdots+\sigma_{n,k_s})}{2^n} \neq 0.
\end{equation}
\end{remark}

Consider now the coefficients $$c_j(k)=\frac{1}{2^r}\sum_{m=0}^{2^r-1}(-1)^{\binom{m}{k}}e^{\frac{-\pi \sqrt{-1}m j}{2^{r-1}}}$$ with $j>0$.  From Theorem \ref{invariant} we know each $c_j(k)$ is the coefficient of $(1+\zeta)^n$ where $1+\zeta$ is a root of $\Phi_{2^{t+1}}(x-1)$ for some $t=1,2,\cdots,2^{r-1}$.
%We will show that $c_j(k)=0$ if and only if it is the coefficient of $(1+\zeta)^n$ where $1+\zeta$ is a root of $\Phi_{2^{b_i+1}}(x-1)$.
\begin{lemma}
\label{cjis0}
Let $k\geq2$ be an integer with $2$-adic expansion $k=2^{a_1}+\cdots+2^{a_l}$.  Then $c_j(k)=0$ if an only if it is the coefficient of $(1+\zeta)^n$, where $1+\zeta$ is a root of $\Phi_{2^{b+1}}(x-1)$ and $b\neq a_i$ for all $i=1,\cdots,l$, i.e. $2^b$ does not appear in the $2$-adic expansion of $k$.
\end{lemma}

\begin{proof}
Recall that to show that the coefficients of the roots of $\Phi_{2^{t+1}}(x-1)$ are zero is equivalent to show that
\begin{equation}
\label{coeffszeros}
\sum_{m=0}^{2^r-1}(-1)^{\binom{m}{k}}e^{\frac{\pi\sqrt{-1}m}{2^{t}}}=0.
\end{equation}
If $\{2^{b_1},\cdots,2^{b_f}\}=\{1,2,2^2,\cdots,2^{r-1}\}\backslash\{2^{a_1},\cdots,2^{a_l}\}$, then equation (\ref{negativeM}) implies,
\begin{eqnarray}
\label{equivsum}
\sum_{m=0}^{2^r-1}(-1)^{\binom{m}{k}}e^{\frac{\pi\sqrt{-1}m}{2^{t}}}&=&-2\sum_{(\delta_1,\cdots,\delta_f)\in \mathbb{F}_2^f}e^{\frac{\pi\sqrt{-1}}{2^{t}}(k+\delta_1 2^{b_1}+\cdots+\delta_t 2^{b_f})}\\ \nonumber
&=&-2e^{\frac{\pi\sqrt{-1}k}{2^{t}}}\sum_{(\delta_1,\cdots,\delta_f)\in \mathbb{F}_2^f}e^{\frac{\pi\sqrt{-1}}{2^{t}}(\delta_1 2^{b_1}+\cdots+\delta_t 2^{b_f})}.
\end{eqnarray}
Thus, (\ref{coeffszeros}) holds if and only if $e^{\frac{\pi\sqrt{-1}}{2^t}}$ is a root of 
\begin{equation}
\label{nicepoly}
\sum_{(\delta_1,\cdots,\delta_f)\in \mathbb{F}_2^f}x^{\delta_1 2^{b_1}+\cdots+\delta_t 2^{b_f}}.
\end{equation}

Consider first $t=b_1$.  If we set $\delta_1=0$ in the last sum of (\ref{equivsum}), then we have
\begin{equation}
\label{pos}
\sum_{(\delta_2,\cdots,\delta_f)\in \mathbb{F}_2^{f-1}}e^{\frac{\pi\sqrt{-1}}{2^{b_1}}(\delta_2 2^{b_2}+\cdots+\delta_t 2^{b_f})}.
\end{equation}
However, if we set $\delta_1=1$, then we have
\begin{equation}
\label{neg}
-\sum_{(\delta_2,\cdots,\delta_f)\in \mathbb{F}_2^{f-1}}e^{\frac{\pi\sqrt{-1}}{2^{b_1}}(\delta_2 2^{b_2}+\cdots+\delta_t 2^{b_f})}.
\end{equation}
We conclude that (\ref{coeffszeros}) holds for $t=b_1$, i.e. the $2^{b_1}$ coefficients related to the roots of $\Phi_{2^{b_1}+1}(x-1)$ are zero.   
Repeat this argument with $t=b_2, \cdots, b_f$ to conclude that the coefficients related to the roots of $\Phi_{2^{b_i+1}}(x-1)$, $i=1,\cdots,f$ are zero.  Since (\ref{nicepoly}) is of degree $d=2^{b_1}+\cdots+2^{b_f}$, then only $d$ of the coefficients $c_j(k)$ can be zero.  Since we already found $d$ coefficients that are zero, then we conclude that these are all of them.  The claim follows.
\end{proof}

Lemmas \ref{c0is0} and \ref{cjis0} are put together in the following theorem.  The function $\epsilon(n)$ used in the theorem is defined as
\begin{equation}
\epsilon(n)=\left\{\begin{array}{cl}
    0, & \text{if }n \text{ is a power of 2,}  \\
    1, & \text{otherwise.}
       \end{array}\right.
\end{equation}
%\begin{remark}
%The sequence $\{\epsilon(n)\}_{n\in \mathbb{N}}$ is known as the Fredholm-Rueppel sequence (sequence $A036987$ in Sloane).
%\end{remark}

\begin{theorem}
\label{charpolyn}
Let $k$ be a natural number and $P_k(x)$ be the characteristic polynomial associated to the minimal linear recurrence with integer coefficients that $\{S(\sigma_{n,k})\}_{n\in\mathbb{N}}$ satisfies.  Let $\bar{k}=2\lfloor k/2\rfloor+1$.  We know $\bar{k}$ has a $2$-adic expansion of the form
\begin{equation}
\bar{k}=1+2^{a_1}+2^{a_2}+\cdots+2^{a_l},
\end{equation}
where the last exponent is given by $a_l=\lfloor \log_2(\bar{k})\rfloor.$  Then $P_k(x)$ equals
\begin{equation}
(x-2)^{\epsilon(k)}\prod_{j=1}^l \Phi_{2^{a_j+1}}(x-1).
\end{equation}
In particular, the degree of the minimal linear recurrence that $\{S(\sigma_{n,k})\}_{n\in\mathbb{N}}$ satisfies is equal to $2\lfloor k/2 \rfloor + \epsilon(k)$.
\end{theorem}

\noindent
Theorem \ref{charpolyn} can be generalized to the case $\{S(\sigma_{n,k_1}+\cdots+\sigma_{n,k_s})\}_{n\in \mathbb{N}}$.  Define the ``OR" operator $\vee$ on $\mathbb{F}_2$ as
\begin{eqnarray}
0 \vee 0 &=&0 \\ \nonumber
0 \vee 1 &=&1 \\ \nonumber
1 \vee 0 &=&1 \\ \nonumber
1 \vee 1 &=&1.
\end{eqnarray}
Extend $\vee$ to $\mathbb{N}$ by letting $m\vee n$ be the natural number obtained by applying $\vee$ coordinatewise to the binary digits of $n$ and $m$.  For example,
\begin{eqnarray}
4\vee 6 &=& (0\cdot 1+0\cdot2+1\cdot 2^2) \vee (0\cdot 1+1\cdot 2+1\cdot 2^2)\\ \nonumber
 &=& (0\vee0)\cdot 1+(0\vee 1)\cdot 2+(1\vee1)\cdot 2^2 =6.
\end{eqnarray}
and
\begin{eqnarray}
3\vee 8 &=& (1\cdot 1+1\cdot2+0\cdot 2^2+0\cdot2^3) \vee (0\cdot 1+0\cdot 2+0\cdot 2^2+1\cdot 2^3)\\ \nonumber
 &=& (1\vee0)\cdot 1+(0\vee 1)\cdot 2+(0\vee0)\cdot 2^2+(0\vee1)\cdot 2^3 =11.
\end{eqnarray}
Next is a generalization of Theorem \ref{charpolyn}.
\begin{theorem}
\label{charpolyvariousk}
Let $1\leq k_1 <k_2<\cdots<k_s$ be fixed integers and $P_{k_1,\cdots,k_s}(x)$ be the characteristic polynomial associated to the minimal linear recurrence with integer coefficients that $\{S(\sigma_{n,k_1}+\cdots+\sigma_{n,k_s})\}_{n\in\mathbb{N}}$ satisfies. Let $\bar{k}= 2\lfloor (k_1\vee\cdots\vee k_s)/2\rfloor+1$.  We know $\bar{k}$ has a $2$-adic expansion of the form
\begin{equation}
\bar{k}=1+2^{a_1}+2^{a_2}+\cdots+2^{a_l},
\end{equation}
where the last exponent is given by $a_l=\lfloor \log_2(\bar{k})\rfloor.$  Then $P_{k_1,\cdots,k_s}(x)$ divides the polynomial
\begin{equation}
(x-2)\prod_{j=1}^l \Phi_{2^{a_j+1}}(x-1).
\end{equation}
\end{theorem}

\begin{proof}
The proof is similar to the one of Lemma \ref{cjis0}.  Let $\bar{k}=2\lfloor (k_1\vee\cdots\vee k_s)/2\rfloor+1$ and $r=\lfloor\log_2(\bar{k})\rfloor+1$.  Define 
\begin{equation}
N=\left\{m \in \{1,2,3,\cdots,2^r-1\}: \binom{m}{k_1}+\cdots+\binom{m}{k_s}\text{ is odd}\right\}.
\end{equation}
Suppose $2^b \in \{2,2^2,\cdots,2^{r-1}\}$ is such that $2^b$ does not appear in the $2$-adic expansion  of $\bar{k}$.  We will show that
\begin{equation}
\label{cond0}
\sum_{m=0}^{2^r-1}(-1)^{\binom{m}{k_1}+\cdots+\binom{m}{k_s}}e^{\frac{\pi\sqrt{-1}}{2^b}}=0,
\end{equation}
which implies that the coefficients related to the roots of $x^{2^b}+1$ are all zero. 

Observe that
\begin{eqnarray}
\sum_{m=0}^{2^r-1}(-1)^{\binom{m}{k_1}+\cdots+\binom{m}{k_s}}e^{\frac{\pi\sqrt{-1}}{2^b}}&=&-2\sum_{m\in N} e^{\frac{\pi\sqrt{-1}m}{2^{b}}}.
\end{eqnarray}
Suppose $m \in N$ is such that $2^b$ does not appear in the $2$-adic expansion of $m$.  Note that equation (\ref{negativeM}) implies $m+2^b \in N$. Thus, the same argument as in (\ref{pos}) and (\ref{neg}) imply that (\ref{cond0}) is true.  Hence, the claim follows.
\end{proof}

The following example presents a case when Theorem \ref{charpolyvariousk} is tight.
\begin{example}
\label{extight}
Consider $k_1=6$ and $k_2=17$.  Note that $2\lfloor(6\vee17)/2\rfloor+1=23=1+2+4+16$.  In this case, the characteristic polynomial associated to $\{S(\sigma_{n,6}+\sigma_{n,17})\}_{n\in \mathbb{N}}$ is $P_{6,17}(x)=(x-2)\Phi_4(x-1)\Phi_8(x-1)\Phi_{32}(x-1)$.  This is the best case scenario of Theorem \ref{charpolyvariousk}, i.e we have equality rather than just divisibility.  Also, note that in this case the recurrence given by Theorem \ref{invariant} is of degree $31$, while the minimal linear recurrence is of degree $23$.
\end{example}
The next example presents a case in which Theorem \ref{charpolyvariousk} improves the degree of the homogeneous linear recurrence provided by Theorem \ref{invariant}.  However it is not the best possible degree.
\begin{example}
Consider $k_1=3$, $k_2=5$, and $k_3=17$.  We have $2\lfloor(3\vee5\vee17)/2\rfloor+1=23$. In this case, the characteristic polynomial of the minimal recurrence is $P_{3,5,17}(x)=(x-2)\Phi_{32}(x-1)$.  It divides $(x-2)\Phi_4(x-1)\Phi_8(x-1)\Phi_{32}(x-1)$ as Theorem \ref{charpolyvariousk} predicted, but are clearly not equal. 
The factors $\Phi_4(x-1)$ and $\Phi_8(x-1)$ do not appear in $P_{3,5,17}(x)$.  This means that the coefficients $c_j(3,5,17)$ related to the roots of $\Phi_4(x)=x^2+1$ and $\Phi_8(x)=x^4+1$ are zero.  However, since $2$ and $4$ appear in the $2$-adic expansion of $23$, then Theorem \ref{charpolyvariousk} cannot detect this.
\end{example}

We now provide bounds on the degree of the minimal linear recurrence that $\{S(\sigma_{n,k_1}+\cdots+\sigma_{n,k_s})\}_{n\in \mathbb{N}}$ satisfies.  We start with the following theorem.

\begin{theorem}
\label{alwaysthere}
Suppose $1\leq k_1<\cdots<k_s$ are integers.  Let $r=\lfloor \log_2(k_s) \rfloor+1$.  Then $\Phi_{2^r}(x-1)$ divides $P_{k_1,\cdots,k_s}(x)$, the characteristic polynomial associated to $\{S(\sigma_{n,k_1}+\cdots+\sigma_{n,k_s})\}_{n\in \mathbb{N}}$.
\end{theorem}

\begin{proof}
Note that the theorem will follow if we show that
\begin{equation}
\sum_{m=0}^{2^r-1}(-1)^{\binom{m}{k_1}+\cdots+\binom{m}{k_s}}e^{\frac{\pi\sqrt{-1} m}{2^{r-1}}}\neq 0.
\end{equation}
This is equivalent to show that $x^{2^{r-1}}+1$ does not divide
\begin{equation}
\label{poly}
\sum_{m=0}^{2^r-1}(-1)^{\binom{m}{k_1}+\cdots+\binom{m}{k_s}}x^m.
\end{equation} 

We present the core of our proof with a particular example. The general case will follow in a similar manner.
Consider the case $k_1=3$, $k_2=5$, and $k_3=10$.  Then (\ref{poly}) equals,
\begin{equation}
\label{expoly}
1 + x + x^2 - x^3 + x^4 - x^5 + x^6 + x^7 + x^8 + x^9 - x^{10} + x^{11} + \
x^{12} - x^{13} - x^{14} - x^{15}.
\end{equation}
Look at the sign of $x^j$ for $j=8,9,\cdots, 15$.  If $x^j$ and $x^{j-8}$ have the same sign, then leave the sign of $x^j$ as it is.  Otherwise, change the sign of $x^j$.  After doing this, we get
\begin{equation}
1 + x + x^2 - x^3 + x^4 - x^5 + x^6 + x^7 + x^8 + x^9 + x^{10} - x^{11} + \
x^{12} - x^{13} + x^{14} + x^{15},
\end{equation}
which equals 
\begin{equation}
\label{divisible}
(1+x^8)(1 + x + x^2 - x^3 + x^4 - x^5 + x^6 + x^7).
\end{equation}
Of course, in order to get (\ref{expoly}) back, we need to add to (\ref{divisible}) two times the terms for which we changed their signs:
\begin{equation}
(1+x^8)(1 + x + x^2 - x^3 + x^4 - x^5 + x^6 + x^7)-2x^{10}+2x^{11}-2x^{14}-2x^{15}.
\end{equation}
This last polynomial equals
\begin{equation}
(1+x^8)(1 + x + x^2 - x^3 + x^4 - x^5 + x^6 + x^7)+2x^8(-x^{2}+x^{3}-x^{6}-x^{7}).
\end{equation}
In general,
\begin{eqnarray}
\sum_{m=0}^{2^r-1}(-1)^{\binom{m}{k_1}+\cdots+\binom{m}{k_s}}x^m &=&(x^{2^{r-1}}+1)\left(\sum_{m=0}^{2^{r-1}-1}(-1)^{\binom{m}{k_1}+\cdots+\binom{m}{k_s}}x^m\right)\\ \nonumber
& &+2x^{2^{r-1}}q(x),
\end{eqnarray}
where $q(x)$ is a polynomial of degree at most $2^{r-1}-1$.  We conclude that $x^{2^{r-1}}+1$ does not divide (\ref{poly}) and so the claim follows.
\end{proof}

\begin{corollary}
\label{boundsdegree}
Let $1\leq k_1<\cdots<k_s$ be integers.  Let $D$ be the degree of the minimal homogeneous linear recurrence with integer coefficients that $\{S(\sigma_{n,k_1}+\cdots+\sigma_{n,k_s})\}_{n\in \mathbb{N}}$ satisfies.  Then $2^{\lfloor\log_2(k_s)\rfloor}\leq D \leq 2\lfloor (k_1 \vee\cdots\vee k_s)/2\rfloor+1$.
\end{corollary}

\begin{proof}
Note that the upper bound follows from Theorem \ref{charpolyvariousk} while the lower bound is a consequence of Theorem \ref{alwaysthere}.
\end{proof}
Note that Corollary \ref{boundsdegree} is an improvement of Theorem \ref{invariant} with respect to the degree $D$ of the minimal homogeneous linear recurrence with integer coefficients that $\{S(\sigma_{n,k_1}+\cdots+\sigma_{n,k_s})\}_{n\in \mathbb{N}}$ satisfies.  From Theorem \ref{invariant} we can only infer that $D\leq 2^r-1$, where $r=\lfloor \log_2(k_s) \rfloor +1$. However, now we know that $2^{\lfloor\log_2(k_s)\rfloor}\leq D \leq 2\lfloor (k_1 \vee\cdots\vee k_s)/2\rfloor+1$ and $2\lfloor (k_1 \vee\cdots\vee k_s)/2\rfloor+1\leq 2^r-1$.  Also, example \ref{extight} shows that the upper bound of Corollary \ref{boundsdegree} can be attained.  In the next theorem, we show that when $k_s$ (the highest degree) is a power of two, then the lower bound is tight.     
\begin{theorem}
Suppose $1\leq k_1<k_2<\cdots<k_s$ are fixed integers with $k_s=2^{r-1}$ a power of two.  Let $P_{k_1,k_2,\cdots,2^{r-1}}(x)$ be the characteristic polynomial associated to the minimal linear recurrence that $\{S(\sigma_{n,k_1}+\sigma_{n,k_2}+\cdots+\sigma_{n,2^{r-1}})\}_{n\in \mathbb{N}}$ satisfies.  Then
\begin{equation}
P_{k_1,k_2,\cdots,2^{r-1}}(x)=\Phi_{2^r}(x-1)=2+\sum _{m=1}^{2^{r-1}} (-1)^m \binom{2^{r-1}}{m} x^m.
\end{equation}
In particular, $\deg(P_{k_1,k_2,\cdots,2^{r-1}}(x))=2^{r-1}=2^{\lfloor \log_2(k_s)\rfloor}.$
\end{theorem}

\begin{proof}
The theorem will follow if we show that $c_0(k_1,k_2,\cdots,2^{r-1})=0$, and
\begin{equation}
\label{almostallzero}
\sum_{m=0}^{2^r-1}(-1)^{\binom{m}{k_1}+\binom{m}{k_2}+\cdots+\binom{m}{2^{r-1}}}e^\frac{\pi \sqrt{-1} m}{2^j}=0,
\end{equation}
for each $j=1,2,\cdots,r-2$, and
\begin{equation}
\label{coeffnotzero}
\sum_{m=0}^{2^r-1}(-1)^{\binom{m}{k_1}+\binom{m}{k_2}+\cdots+\binom{m}{2^{r-1}}}e^\frac{\pi \sqrt{-1} m}{2^{r-1}}\neq 0.
\end{equation}

From Theorem \ref{alwaysthere} we know that (\ref{coeffnotzero}) holds true.  Now, the coefficient $c_0(k_1,\cdots,k_{s-1},2^{r-1})$ is zero if and only if 
\begin{equation}
\label{c0zero}
\sum_{m=0}^{2^r-1}(-1)^{\binom{m}{k_1}+\cdots+\binom{m}{k_{s-1}}+\binom{m}{2^{r-1}}}=0.
\end{equation}
From (\ref{frombino}) we see that the period of $(-1)^{\binom{m}{k_1}+\cdots+\binom{m}{k_{s-1}}}$ is a proper divisor of $2^r$, so
\begin{equation}
\sum_{m=0}^{2^{r-1}-1}(-1)^{\binom{m}{k_1}+\cdots+\binom{m}{k_{s-1}}}=\sum_{m=2^{r-1}}^{2^r-1}(-1)^{\binom{m}{k_1}+\cdots+\binom{m}{k_{s-1}}}.
\end{equation}
However, 
\begin{equation}
\label{pm}
(-1)^{\binom{m}{2^{r-1}}}=\left\{
\begin{array}{rl}
 1, & \text{if } m\leq 2^{r-1}-1 \\
 -1, & \text{if }m\geq 2^{r-1}.
\end{array}
\right.
\end{equation}
Thus, (\ref{c0zero}) holds and so $c_0(k_1,\cdots,k_{s-1},2^{r-1})=0$.  Similarly, the periodicity of $(-1)^{\binom{m}{k_1}+\cdots+\binom{m}{k_{s-1}}}$ and $e^{\frac{\pi\sqrt{-1}m}{2^j}}$ implies
\begin{equation}
\label{equalsums}
\sum_{m=0}^{2^{r-1}-1}(-1)^{\binom{m}{k_1}+\cdots+\binom{m}{k_{s-1}}}e^{\frac{\pi\sqrt{-1}m}{2^j}}=\sum_{m=2^{r-1}}^{2^r-1}(-1)^{\binom{m}{k_1}+\cdots+\binom{m}{k_{s-1}}}e^{\frac{\pi\sqrt{-1}m}{2^j}}.
\end{equation}
So, (\ref{pm}) and (\ref{equalsums}) imply (\ref{almostallzero}). This concludes the proof.
\end{proof}

We conclude this section with the following result, which shows that when $k_s$ is a power of two, then, as $n$ increases, $|S(\sigma_{n,k_1}+\cdot+\sigma_{n,k_s})|$ is much smaller than $2^n$. 

\begin{corollary}
Suppose $1\leq k_1<k_2<\cdots<k_s$ are fixed integers with $k_s=2^{r-1}$ a power of two.  Then, for $0\leq j \leq 2^r-1$, $c_j(k_1,\cdots, k_{s-1},2^{r-1})\neq 0$
if and only if $j$ is odd.  In particular, $c_0(k_1,\cdots, k_{s-1},2^{r-1})= 0$.
\end{corollary}
%%%%%%%%%%%%%%%%%%%%%%%%%%%%%%%%%%%%%%%%%%%%%%%%%%%%%%%%%%%
% Section: Asymptotic behavior
%%%%%%%%%%%%%%%%%%%%%%%%%%%%%%%%%%%%%%%%%%%%%%%%%%%%%%%%%%%
\section{Asymptotic behavior}

In this section we discuss the asymptotic behavior of $\{S(\sigma_{n,k_1}+\cdots+\sigma_{n,k_s})\}_{n\in \mathbb{N}}$.  Note that Theorem \ref{invariant} implies that
\begin{equation}
\lim_{n\to \infty}\frac{S(\sigma_{n,k_1}+\cdots+\sigma_{n,k_s})}{2^n}=c_0(k_1,\cdots, k_s)=\frac{1}{2^r}\sum_{m=0}^{2^r-1}(-1)^{\binom{m}{k_1}+\cdots+\binom{m}{k_s}}.
\end{equation}
Thus, we study $c_0(k_1,\cdots, k_s)$ first.  

We already discussed the case of one elementary symmetric polynomial $\{S(\sigma_{n,k})\}_{n\in \mathbb{N}}$, see (\ref{valc0}): 
\begin{equation}
\label{formonek}
c_0(k)=\frac{2^{w_2(k)-1}-1}{2^{w_2(k)}}.
\end{equation}
For instance, we know that $c_0(k)\geq 0$, and the equality holds if and only if $k$ is a power of two. 

The method of inclusion-exclusion can be used to get a formula in the case that we have more than one symmetric polynomial.  For example, in the case of two elementary symmetric polynomials $\{S(\sigma_{n,k_1}+\sigma_{n,k_2})\}_{n\in \mathbb{N}}$, we have
\begin{eqnarray}
\label{twoks}
c_0(k_1,k_2) &=& 1-2^{1-w_2(k_1)}-2^{1-w_2(k_2)}+2^{2-w_2(k_1 \vee k_2)}
\end{eqnarray}
The reader can check that this formula implies $c_0(k_1,k_2)\geq 0$, with equality if and only if $w_2(k_1 \vee k_2)=w(k_1)+w_2(k_2)$ and $w_2(k_i)=1$, where $i=1$ or $i=2$.  We start the general case with the following lemma.
\begin{lemma}
\label{setofnegatives}
Suppose that $1\leq k_1<k_2<\cdots<k_s$ are integers.  Define
\begin{equation}
N(k_{i_1},\cdots k_{i_j})=\left\{m \in \{1,2,3,\cdots,2^r-1\}: \binom{m}{k_{i_1}}+\cdots+\binom{m}{k_{i_j}} \text{ is odd} \right\}.
\end{equation}
Then,
\begin{eqnarray}\nonumber
\#N(k_1,k_2,\cdots, k_s)&=&\sum_{i=1}^s \#N(k_i) -2\sum_{i_1<i_2} \#\left(N(k_{i_1})\cap N(k_{i_2})\right)\\
& &+4\sum_{i_1<i_2<i_3} \#\left(N(k_{i_1})\cap N(k_{i_2})\cap N(k_{i_2})\right)-\cdots\\ \nonumber
& &+(-1)^{s-1}2^{s-1}\#\left(N(k_1)\cap N(k_2)\cap\cdots\cap N(k_s)\right).
\end{eqnarray}
\end{lemma}
\begin{proof}
Note that
\begin{equation}
\binom{m}{k_1}+\cdots+\binom{m}{k_s}
\end{equation}
is odd exactly when the amount of odd summands is itself an odd number (trivial).  In terms of our sets $N(k_i)$, this implies that $N(k_1,\cdots, k_s)$ is obtained by including all the intersections of an odd amount of the sets $N(k_i)$, while excluding all the intersections of an even amount of them. For example, the case of four $k$'s can be represented by the Venn diagram in Figure \ref{venn4}.  In this case, we want to include the shaded regions and exclude the white ones.  
\begin{figure}[h!]
\caption{Representation of the case of four $k$'s.}
\label{venn4}
\includegraphics[width=1.35in]{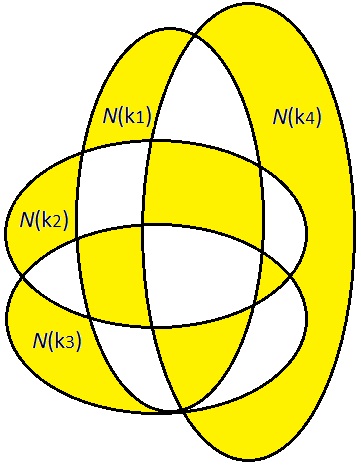}
\end{figure}

We start by adding $\#N(k_1)+\#N(k_2)+\cdots+\#N(k_s)$. Then, we proceed to take out all the intersections of two sets $N(k_i)\cap N(k_j)$, $i\neq j$.  In this case, each of them has been added twice in our previous sum.  Therefore, to take them out, we should add $-2\#(N(k_1)\cap N(k_2))-2\#(N(k_1)\cap N(k_3))-\cdots-2\#(N(k_{s-1})\cap N(k_s))$ to the previous sum.  So, now we have 
\begin{equation}
\label{taketwoint}
\sum_{i=1}^s\#N(k_i)-2\sum_{i_1<i_2}\#(N(k_{i_1})\cap N(k_{i_2})).
\end{equation}
This takes care of all intersections of two sets.  Now, we need to add all intersections of three sets $\#(N(k_{i_1})\cap N(k_{i_2})\cap N(k_{i_3}))$.  Each of them have been added three times by the first sum and subtracted six times by the second sum.  Thus, in order to add them into the equation, we have to add each of them four times to (\ref{taketwoint}).  Doing this, we now have
\begin{eqnarray}
\sum_{i=1}^s\#N(k_i)&-&2\sum_{i_1<i_2}\#(N(k_{i_1})\cap N(k_{i_2}))\\ \nonumber
&+&4\sum_{i_1<i_2<i_3}\#(N(k_{i_1})\cap N(k_{i_2})\cap N(k_{i_3})).
\end{eqnarray}
Continue in this matter and use the identity
\begin{equation}
%\label{pm}
\sum_{i=1}^{j-1}(-1)^{i-1}2^{i-1}\binom{j}{i}=\left\{
\begin{array}{ll}
 2^{j-1}, & \text{if {\it j} is even}\\
 -2^{j-1}+1, & \text{if {\it j} is odd,}
\end{array}
\right.
\end{equation}
to get the result.
\end{proof}

Now that we have the above lemma, we are ready to state our general formula for $c_0(k_1,\cdots, k_s)$.
\begin{theorem}
\label{generalc0}
Suppose that $1\leq k_1<k_2<\cdots<k_s$ are integers.  Then 
\begin{eqnarray}\nonumber
c_0(k_1,\cdots,k_s)&=& 1-\sum_{i=1}^s 2^{1-w_2(k_i)}+\sum_{i_1<i_2} 2^{2-w_2(k_{i_1}\vee k_{i_2})}\\
& &-\sum_{i_1<i_2<i_3} 2^{3-w_2(k_{i_1}\vee k_{i_2}\vee k_{i_3})}+\cdots+(-1)^s 2^{s-w_2(k_1 \vee k_2 \vee \cdots \vee k_s)}.
\end{eqnarray}
\end{theorem}
\begin{proof}
Let $r=\lfloor\log_2(k_s)\rfloor+1$. Note that 
\begin{equation}
\label{easy}
c_0(k_1,\cdots,k_s)=\frac{2^r-\#N(k_1,\cdots,k_s)}{2^r}.
\end{equation}

Consider the integer $k_i$.  Suppose its $2$-adic expansion is $k_i=2^{a_1}+2^{a_2}+\cdots+2^{a_l}$. Then, by Corollary \ref{modds} we know that $0\leq m\leq 2^r-1$ is such that $\binom{m}{k_i}$ is odd precisely when
\begin{equation}
\label{negativeMagain}
m=k_i+\delta_1 2^{b_1}+\delta_2 2^{b_2}+\cdots+\delta_t 2^{b_t},
\end{equation}
where $\{2^{b_1},2^{b_2},\cdots, 2^{b_t}\}=\{1,2,2^2,\cdots, 2^{r-1}\} \backslash\{2^{a_1},2^{a_2},\cdots, 2^{a_l}\}$ and $\delta_j$ is either 0 or 1.  This implies that \begin{equation}
\#N(k_i)=2^{r-w_2(k_i)}.
\end{equation}  
Also, (\ref{negativeMagain}) implies that if $i\neq j$, then 
\begin{equation}
\#(N(k_i)\cap N(k_j))=2^{r-w_2(k_i\vee k_j)},
\end{equation} 
and in general, if $1\leq i_1<i_2<\cdots i_t\leq s$, then 
\begin{equation}
\#(N(k_{i_1})\cap N(k_{i_2}) \cap \cdots \cap N(k_t))=2^{r-w_2(k_{i_1}\vee k_{i_2}\vee\cdots\vee k_{i_t})}.
\end{equation}
This, plus equation (\ref{easy}) and Lemma \ref{setofnegatives}, imply the theorem.
\end{proof}

\begin{example}
Suppose $k_1=7$, $k_2=9$, $k_3=2^{10^5}+2^{10^4}$, and $k_4=2^{10^6}+5$, then $c_0(k_1,k_2,k_3,k_4)=1/4$.  In other words,
\begin{equation}
\lim_{n\to \infty}\frac{S(\sigma_{n,k_1}+\sigma_{n,k_2}+\sigma_{n,k_3}+\sigma_{n,k_4})}{2^n}=\frac{1}{4}.
\end{equation}
\end{example}

\begin{example}
Suppose $k_1=31$, $k_2=2^{10^4}+64$ and $k_3=2^{10^4}+32+128$, then $c_0(k_1,k_2,k_3)=45/128$.
\end{example}

We now use the above theorem to provide a family of symmetric polynomials that are not balanced for sufficiently large $n$. Suppose that $k_1$ and $k_2$ are two positive integers.  We say that $k_1 \preceq k_2$ if each power of two appearing in the $2$-adic expansion of $k_1$ also appears in the $2$-adic expansion of $k_2$.  For example, $10 \preceq 14$ because $10=2+8$ and $14=2+4+8$.

\begin{corollary}
\label{ksinsideoneanother}
Suppose that $k_1 \preceq k_2 \preceq \cdots \preceq k_s$ are positive integers.  Then,
\begin{equation}
c_0(k_1,\cdots,k_s)=1-\Delta(s) 2^{1-w_2(k_s)}-\sum_{j=1}^{\lfloor s/2\rfloor}(2^{w_2(k_{2j}-k_{2j-1})}-1)2^{1-w_2(k_{2j})}.
\end{equation}
Here $\Delta(s)$ equals $0$ if $s$ is even and $1$ otherwise; in other words, $\Delta(s)= s \mod 2$.
\end{corollary}
\begin{proof}
This follows directly from Theorem \ref{generalc0} and the equality 
$$w_2(k_i)=w_2(k_i-k_{i-1})+w_2(k_{i-1}). $$
\end{proof}

\begin{theorem}
Suppose that $k_1 \preceq k_2 \preceq \cdots \preceq k_s$ are positive integers.  Then,
\begin{equation}
c_0(k_1,\cdots,k_s)> 0.
\end{equation}
In particular, $\{S(\sigma_{n,k_1}+\cdots+\sigma_{n,k_s}\}_{n \in \mathbb{N}}$ is asymptotically not balanced.
\end{theorem}

\begin{proof}
Corollary \ref{ksinsideoneanother} implies that $c_0(k_1,\cdots,k_s)>0$ if and only if
\begin{equation}
\frac{\Delta(s)}{2^{w_2(k_s)}}+\sum_{j=1}^{\lfloor s/2\rfloor}\frac{2^{w_2(k_{2j}-k_{2j-1})}-1}{2^{w_2(k_{2j})}}<\frac{1}{2}.
\end{equation}
Since $k_1\preceq k_2\preceq\cdots\preceq k_s$, then we have the inequality
\begin{equation}
\label{relki}
w_2(k_i)\geq 1+w_2(k_{i-1})
\end{equation}
and the equality 
\begin{equation}
\label{equalityki}
w_2(k_i)=w_2(k_i-k_{i-1})+w_2(k_{i-1}).  
\end{equation}
Note that (\ref{relki}) and (\ref{equalityki}) imply
\begin{eqnarray}\nonumber
\frac{\Delta(s)}{2^{w_2(k_s)}}+\sum_{j=1}^{\lfloor s/2\rfloor}\frac{2^{w_2(k_{2j}-k_{2j-1})}-1}{2^{w_2(k_{2j})}}&=&\frac{\Delta(s)}{2^{w_2(k_s)}}+\frac{1}{2^{w_2(k_1)}}-\frac{1}{2^{w_2(k_2)}}+\sum_{j=2}^{\lfloor s/2\rfloor}\frac{2^{w_2(k_{2j}-k_{2j-1})}-1}{2^{w_2(k_{2j})}}\\
&\leq&\frac{\Delta(s)}{2^{w_2(k_s)}}+\frac{1}{2^{w_2(k_1)}}-\frac{1}{2^{w_2(k_2)}}+\frac{1}{2^{w_2(k_2)}}\sum_{j=2}^{\lfloor s/2\rfloor}\frac{1}{2^{2j-3}}\\\nonumber
&<&\frac{1}{2^{w_2(k_1)}}-\frac{1}{2^{w_2(k_2)}}+\frac{1}{2^{w_2(k_2)}}\sum_{j=2}^{\infty}\frac{1}{2^{2j-3}}\\ \nonumber
&=&\frac{1}{2^{w_2(k_1)}}+\left(\frac{2}{3}-1\right)\frac{1}{2^{w_2(k_2)}}<\frac{1}{2}.
\end{eqnarray}
This finishes the proof.
\end{proof}

We now turn our attention to the case when $c_0(k_1,\cdots,k_s)=0$, which happens if and only if $2$ is not a root of the characteristic polynomial associated to $\{S(\sigma_{n,k_1}+\cdots+\sigma_{n,k_s}\}_{n \in \mathbb{N}}$.  In this case
\begin{equation}
\lim_{n\to\infty}\frac{S(\sigma_{n,k_1}+\cdots+\sigma_{n,k_s})}{2^n}=0,
\end{equation}
because $|S(\sigma_{n,k_1}+\cdots+\sigma_{n,k_s})|$ is exponentially smaller than $2^n$.  However, aside from the size of $S(\sigma_{n,k_1}+\cdots+\sigma_{n,k_s})$ with respect to $2^n$, knowing that $c_0(k_1,\cdots,k_s)=0$ does not give a real sense of the behavior of $S(\sigma_{n,k_1}+\cdots+\sigma_{n,k_s})$ as $n$ increases.

Now, when $c_0(k_1,\cdots,k_s)=0$, the biggest modulus of the roots of the characteristic polynomial associated to $\{S(\sigma_{n,k_1}+\cdots+\sigma_{n,k_s}\}_{n \in \mathbb{N}}$ is $2\cos(\pi/2^r)$.  This modulus is obtained at the roots $1+e^{\pi\sqrt{-1}/(2^{r-1})}$ and $1+e^{-\pi\sqrt{-1}/(2^{r-1})}$.  Thus, as $n$ increases,
\begin{equation}
\frac{S(\sigma_{n,k_1}+\cdots+\sigma_{n,k_s})}{(2\cos(\pi/2^r))^n}
\end{equation}
approaches
\begin{equation}
\frac{c_1(k_1,\cdots,k_s)(1+e^{\frac{\pi\sqrt{-1}}{2^{r-1}}})^n+c_{2^r-1}(k_1,\cdots,k_s)(1+e^{\-\frac{\pi\sqrt{-1}}{2^{r-1}}})^n}{(2\cos(\pi/2^r))^n}.
\end{equation}
Let $c_i=c_i(k_1,\cdots,k_s)$ and $\xi = 1+e^{\frac{\pi\sqrt{-1}}{2^{r-1}}}$. Note that $c_{2^r-1}=\bar{c_1}$. Since
\begin{equation}
1+e^{\pm\frac{\pi\sqrt{-1}}{2^{r-1}}} = 2\cos\left(\frac{\pi}{2^r}\right)e^{\pm\frac{\pi\sqrt{-1}}{2^r}},
\end{equation}
then
\begin{eqnarray}\nonumber
\label{asymp1}
c_1\xi^n+\bar{c_1}\bar{\xi}^n&=&
\frac{(2\cos(\pi/2^r))^n}{2^r}\sum_{m=0}^{2^r-1}(-1)^{\binom{m}{k_1}+\cdots+\binom{m}{k_s}}\left(e^{\frac{(n-2m)\pi\sqrt{-1}}{2^r}}+e^{-\frac{(n-2m)\pi\sqrt{-1}}{2^r}}\right)\\
&=&\frac{(2\cos(\pi/2^r))^n}{2^{r-1}}\sum_{m=0}^{2^r-1}(-1)^{\binom{m}{k_1}+\cdots+\binom{m}{k_s}}\cos\left(\frac{(n-2m)\pi}{2^r}\right).
\end{eqnarray}
Note that (\ref{asymp1}) is not the zero function, because
\begin{equation}
\sum_{m=0}^{2^r-1}(-1)^{\binom{m}{k_1}+\cdots+\binom{m}{k_s}}\cos\left(\frac{(n-2m)\pi}{2^r}\right)
\end{equation}
is the real part of 
\begin{equation}
\label{totakerealpart}
e^{\frac{\pi\sqrt{-1}n}{2^r}}\sum_{m=0}^{2^r-1}(-1)^{\binom{m}{k_1}+\cdots+\binom{m}{k_s}}e^{-\frac{m\pi\sqrt{-1}}{2^{r-1}}}
\end{equation}
and from the proof of Theorem \ref{boundsdegree} we know that the sum in (\ref{totakerealpart}) is a non-zero constant.  We conclude that
\begin{eqnarray} \nonumber
\label{asympexprc0is0}
S(\sigma_{n,k_1}+\cdots+\sigma_{n,k_s})&=&\frac{(2\cos(\pi/2^r))^n}{2^{r-1}}\sum_{m=0}^{2^r-1}(-1)^{\binom{m}{k_1}+\cdots+\binom{m}{k_s}}\cos\left(\frac{(n-2m)\pi}{2^r}\right)\\
& & +\,O\left(\left(2\cos\left(\frac{\pi}{2^{r-1}}\right)\right)^n\right).\\\nonumber
\end{eqnarray}

\begin{remark}
Note that 
\begin{equation}
\sum_{m=0}^{2^r-1}(-1)^{\binom{m}{k_1}+\cdots+\binom{m}{k_s}}\cos\left(\frac{(n-2m)\pi}{2^r}\right),
\end{equation}
is not identically zero, regardless if $c_0(k_1,\cdots,k_s)$ is zero.  Hence, the asymptotic expansion of $S(\sigma_{n,k_1}+\cdots+\sigma_{n,k_s})$ is 
\begin{eqnarray} \nonumber
\label{asympexpr}
S(\sigma_{n,k_1}+\cdots+\sigma_{n,k_s})&=&c_02^n+\frac{(2\cos(\pi/2^r))^n}{2^{r-1}}\sum_{m=0}^{2^r-1}(-1)^{\binom{m}{k_1}+\cdots+\binom{m}{k_s}}\cos\left(\frac{(n-2m)\pi}{2^r}\right)\\
& & +\,O\left(\left(2\cos\left(\frac{\pi}{2^{r-1}}\right)\right)^n\right).\\\nonumber
\end{eqnarray}
\end{remark}

\begin{remark}
If $c_0(k_1,\cdots,k_s)=0$, then we define the function $\text{Error}_n(k_1,\cdots,k_s)$ as
\begin{eqnarray} \nonumber
\text{Error}_n(k_1,\cdots,k_s)&=&\frac{S(\sigma_{n,k_1}+\cdots+\sigma_{n,k_s})}{(2\cos(\pi/2^r))^n}\\ 
& & -\frac{1}{2^{r-1}}\sum_{m=0}^{2^r-1}(-1)^{\binom{m}{k_1}+\cdots+\binom{m}{k_s}}\cos\left(\frac{(n-2m)\pi}{2^r}\right). 
\end{eqnarray}
\end{remark}
We point out that (\ref{asympexpr}) is a consequence of Theorem \ref{invariant}.  Moreover, it is not hard to continue refining this formula and re-write $S(\sigma_{n,k_1}+\cdots+\sigma_{n,k_s})$ completely in terms of cosine (which will be a restatement of Theorem \ref{invariant}).  However, we stress out that we are proving that
\begin{equation}
\sum_{m=0}^{2^r-1}(-1)^{\binom{m}{k_1}+\cdots+\binom{m}{k_s}}\cos\left(\frac{(n-2m)\pi}{2^r}\right)
\end{equation}
is not identically zero, so the second term of the asymptotic expansion (\ref{asympexpr}) is always present.  By the discussion presented in section \ref{sectiondegree} we know that in many cases some of the $c_i$'s, $i\neq 1,2^r-1$, are zero. Because of this, we write the asymptotic expansion of $S(\sigma_{n,k_1}+\cdots+\sigma_{n,k_s})$ as in (\ref{asympexpr}) and decide not to continue with a refinement of it.
\begin{example}
Consider $k_1=5$, $k_2=9$, and $k_3=12$.  The reader can check that in this case $c_0(5,9,12)=0$.  Therefore, by (\ref{asympexprc0is0}) we know that as $n$ increases, 
\begin{equation}
\label{ex5912}
\frac{S(\sigma_{n,5}+\sigma_{n,9}+\sigma_{n,12})}{(2\cos(\pi/16))^n}
\end{equation}
approaches
\begin{eqnarray}
\label{ex5912asym}
&&\frac{1}{8}\sum_{m=0}^{15} (-1)^{\binom{m}{5}+\binom{m}{9}+\binom{m}{12}}\cos\left(\frac{(n-2m)\pi}{16}\right)=\\ \nonumber
&& \frac{1}{8} \left(\sqrt{2} \left(\sqrt{2+\sqrt{2}}-1\right) \cos \left(\frac{n \pi }{16}\right)+\left(2+\sqrt{2}+\sqrt{2
   \left(2+\sqrt{2}\right)}\right) \sin \left(\frac{n \pi }{16}\right)\right).
\end{eqnarray}
Here, we simplified the expression using {\it Mathematica}.  In Figure \ref{figureex5912} you can see a graphical representation of this.  The blue dots represents (\ref{ex5912}) and in red is the curve given by (\ref{ex5912asym}).  In Table \ref{error5912} you can see the error term.
\begin{figure}[h!]
\caption{Graphical representation of the asymptotic behavior when $k_1=5$, $k_2=9$, and $k_3=12$.}
\includegraphics[width=2.1in]{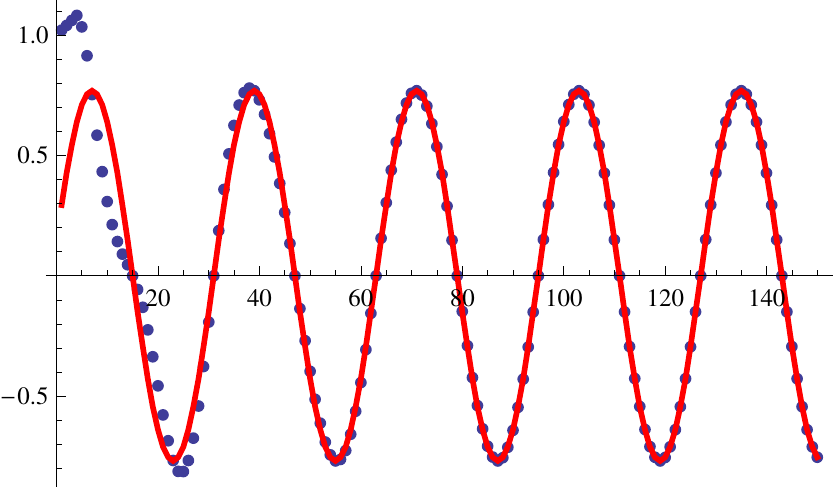}
\label{figureex5912}
\end{figure}
\begin{table}[h!]
\caption{The error term for $k_1=5$, $k_2=9$, and $k_3=12$.}
\label{error5912}
\begin{tabular}{|l|l|}
\hline
 $n$ & $\text{Error}_n(5,9,12)$\\
 \hline
 100 & $0.001530582098$ \\
 200 & $-1.60776038707\times10^{-6}$ \\
 300 & $-9.843230768196\times10^{-9}$ \\
 400 & $1.033957384537\times10^{-11}$ \\
 500 & $6.330222602868\times10^{-14}$ \\
 \hline
\end{tabular}
\end{table}
\end{example}

\begin{example}
Similarly, consider $k_1=2$, $k_2=4$, $k_3=11$, and $k_4=35$.  In this case, we also have $c_0(5,9,12)=0$, so by (\ref{asympexprc0is0}) we know that as $n$ increases, 
\begin{equation}
\label{ex241135}
\frac{S(\sigma_{n,2}+\sigma_{n,4}+\sigma_{n,11}+\sigma_{n,35})}{(2\cos(\pi/64))^n}
\end{equation}
approaches
\begin{equation}
\label{ex241135asym}
\frac{1}{32}\sum_{m=0}^{63} (-1)^{\binom{m}{2}+\binom{m}{4}+\binom{m}{11}+\binom{m}{35}}\cos\left(\frac{(n-2m)\pi}{16}\right).
\end{equation}
In Figure \ref{figureex241135} you can see a graphical representation of this, where the blue dots represents (\ref{ex241135}) and in red is the curve given by (\ref{ex241135asym}).  In Table \ref{error241135} you can see the error term.
\begin{figure}[h!]
\caption{Graphical representation of the asymptotic behavior when $k_1=2$, $k_2=4$, $k_3=11$, and $k_4=35$.}
\includegraphics[width=2.1in]{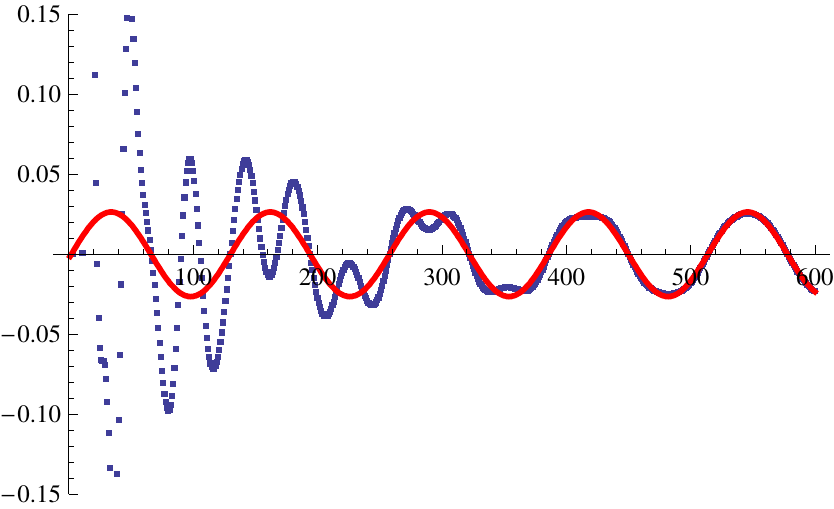}
\label{figureex241135}
\end{figure}
\begin{table}[ht]
\caption{The error term for $k_1=2$, $k_2=4$, $k_3=11$, and $k_3=35$.}
\label{error241135}
\begin{tabular}{|l|l|}
\hline
 $n$ & $\text{Error}_n(2,4,11,35)$\\
 \hline
 $250$ & $-0.014750$ \\
 $500$ & $-0.0012673$ \\
 $750$ & $-0.000024944$ \\
 $1000$ & $7.21779483609288\times10^{-6}$ \\
 $1250$ & $1.01240694303367\times10^{-6}$\\ 
 \hline
\end{tabular}
\end{table}
\end{example}
%%%%%%%%%%%%%%%%%%%%%%%%%%%%%%%%%%%%%%%%%%%%%%%%%%%%%%%%%%%
% Acknowledgments
%%%%%%%%%%%%%%%%%%%%%%%%%%%%%%%%%%%%%%%%%%%%%%%%%%%%%%%%%%%
\bigskip
\noindent
{\bf Acknowledgments.} We would like to thank Claude Carlet for a careful reading of the manuscript and for all his helpful suggestions. \\

%%%%%%%%%%%%%%%%%%%%%%%%%%%%%%%%%%%%%%%%%%%%%%%%%%%%%%%%%%%
% Bibliography
%%%%%%%%%%%%%%%%%%%%%%%%%%%%%%%%%%%%%%%%%%%%%%%%%%%%%%%%%%%
%\bibliography{../AllRef/Bib}
\bibliographystyle{AMSplain}

\end{document}